\documentclass[preprint,12pt]{elsarticle}

\usepackage{amssymb}
\usepackage{amsmath}
\usepackage{amsthm}

\usepackage{graphicx}

\journal{Systems \& Control Letters}

\begin{document}
\newtheorem{thm}{Theorem}
\newtheorem{pro}{Proposition}

\newtheorem{lem}[thm]{Lemma}
\newtheorem{example}{Example}
\newtheorem{defn}{Definition}


\newcommand\scalemath[2]{\scalebox{#1}{\mbox{\ensuremath{\displaystyle #2}}}}

\begin{frontmatter}



\title{Observations on robust diffusive stability and common Lyapunov functions}


\author[1,2]{Blake McGrane-Corrigan} 

\affiliation[1]{organization={School of Biology and Environmental Sciences, 
University College Dublin},
            addressline={Belfield}, 
            city={Dublin 4},
            country={Ireland}}
            
\author[2]{Rafael de Andrade Moral}

\author[2]{Oliver Mason\corref{cor1}} 
\ead{oliver.mason@mu.ie}

\cortext[cor1]{Corresponding Author}
\affiliation[2]{organization={Dept. of Mathematics \& Statistics, Maynooth University},
            city={Maynooth},
            state={Co. Kildare},
            country={Ireland}}

\begin{abstract}
We consider the problem of robust diffusive stability (RDS) for a pair of coupled stable discrete-time positive linear-time invariant (LTI) systems.  We first show that the existence of a common diagonal Lyapunov function is sufficient for RDS and highlight how this condition differs from recent results using linear copositive Lyapunov functions.  We also present an extension of these results, showing that the weaker condition of \emph{joint} linear copositive function existence is also sufficient for RDS.  Finally, we present two results on RDS for extended Leslie matrices arising in population dynamics.  
\end{abstract}




\end{frontmatter}



\section{Introduction}
\label{sec1}
In recent decades, numerous authors have studied problems related to the stability of positive systems.  Such systems, where the state remains nonnegative for all time when starting from nonnegative initial conditions, arise in connection with applications in areas such as population dynamics, economics, and communications \cite{FarRin}. The stability theory of positive systems has several distinct features.  In particular, positive systems often possess stronger stability properties such as $D$-stability or delay-independent stability, and these can be established using classes of Lyapunov functions such as diagonal or linear copositive Lyapunov functions \cite{Mas, AlMas}.  The basic stability properties of positive linear time-invariant (LTI) systems are well understood, and much of the recent work in the field has focused on extensions to systems subject to nonlinearity, time-delay, uncertainty, or time-variation/switching \cite{AlMas, ForVal}.  In the context of switched systems, it is well-known that the existence of a common linear copositive Lyapunov function (CLCLF) or a common diagonal Lyapunov function (CDLF) is sufficient for stability under arbitrary switching \cite{ForVal,MasSho}.  For this reason, many researchers have studied these and related Lyapunov functions.  So-called \emph{joint linear copositive Lyapunov functions} \cite{WuSun} have also been studied, as these give a potentially less conservative condition than CLCLFs.   

Inspired by the recent papers \cite{deLee,Guiv} addressing the impact of movement on the stability of equilibria in population dynamics, we consider the problem of \emph{robust diffusive stability} for a pair of coupled, discrete-time stable positive LTI systems $\Sigma_A: x(t+1)=Ax(t)$, $\Sigma_B: x(t+1) = B x(t)$.  In \cite{Guiv}, the authors focussed on what they term \emph{dispersal-driven growth}, which arises when two stable LTI systems can be destabilised by coupling.  The results in \cite{Guiv} give necessary conditions for dispersal-driven growth based on the \emph{non-existence} of common linear copositive Lyapunov functions (CLCLF) \cite{MasSho}.  They do this by showing that the existence of a CLCLF is a sufficient condition for \emph{robust diffusive stability} (RDS): informally, RDS means that it is impossible to destabilise by coupling.   An extension of this fact to the setting of continuous-time coupled LTI systems with an invariant cone was developed in \cite{deLee}.  Both of these papers describe a novel and intriguing application of  CLCLFs to the stability of diffusively coupled systems. We develop on this theme by considering other types of common Lyapunov functions, which were introduced for the analysis of \emph{switched positive systems}, and investigate how these relate to the RDS property for coupled systems.

The layout of the paper is as follows.  After recalling relevant background results and formally defining the problem in Section \ref{sec:prel}, we prove in Section \ref{sec:2} that CDLF existence gives an alternative sufficient condition for robust diffusive stability when the coupling matrices are diagonal.  In Section \ref{sec:jlf}, we prove that a joint linear copositive Lyapunov function (JLCLF) is sufficient for RDS under a mild irreducibility assumption, thereby giving a less restrictive condition than that given by CLCLF existence.   In Section \ref{sec:Les}, we present results on robust diffusive stability for systems arising in connection with population dynamics where the system matrices $A$, $B$ are extended Leslie matrices.
  
\section{Preliminaries and problem definition}\label{sec:prel}
Our notation is fairly standard. We will use $\mathbb{R}^n_+$ to denote the nonnegative orthant and $\mathbb{R}^{n \times n}_+$ to denote the cone of (entrywise) nonnegative matrices.  For a symmetric matrix $P=P^T$ in $\mathbb{R}^{n \times n}$, $P \succ 0$ ($P \prec 0$) denotes that $P$ is positive definite (negative definite).  For $v \in \mathbb{R}^n$, $v \geq 0$ means that $v_i \geq 0$ for $1 \leq i \leq n$, $v > 0$ means that $v \geq 0$, $v \neq 0$, while $v \gg 0$ means $v_i > 0$ for $1\leq i \leq n$.  For $A \in \mathbb{R}^{n \times n}$, $\rho(A)$ denotes its spectral radius.  $A \in \mathbb{R}^{n \times n}$ is Schur-stable if the spectral radius $\rho(A)$ of $A$ satisfies $\rho(A) < 1$.  

$A$ in $\mathbb{R}^{n \times n}$ is reducible if there exist non-empty, disjoint subsets $I$, $J$ of $\{1, \ldots n\}$ with $I \cup J = \{1,\ldots n\}$ such that $a_{ij} = 0$ for all $i \in I$, $j \in J$.  If $A$ is not reducible, it is \emph{irreducible} \cite{HorJoh}. 

In the next result, we gather several well-known characterisations of Schur-stability for a single nonnegative matrix \cite{FarRin,KasBha}.

\begin{thm}\label{thm:SchStab}
Let $A \in \mathbb{R}^{n \times n}_+$ be a nonnegative matrix.  The following are equivalent:
\begin{enumerate}
\item $A$ is Schur-stable;
\item there exists $v \gg 0$ in $\mathbb{R}^n_+$ with $v^T A \ll v^T$;
\item there is no non-zero $w > 0$ in $\mathbb{R}^n_+$ with $Aw \geq w$;
\item there exists a diagonal matrix $E \succ 0$ with $A^TEA - E \prec 0$;
\item there exists a diagonal matrix $E \succ 0$ with $(A-I)^T E + E(A-I) \prec 0$. 
\end{enumerate}
\end{thm}
The matrix inequality 4 in Theorem \ref{thm:SchStab} is known as the Stein inequality, while inequality 5 is the Lyapunov inequality (for $A-I$).  Note that as $\rho(A) = \rho(A^T)$, another equivalent condition for Schur-stability is the existence of $v \gg 0$ with $A v \ll v$.  

\subsection{Problem definition}

Consider two positive discrete-time LTI systems
\[\Sigma_A: x(t+1) = A x(t), \;\;\; \Sigma_B: y(t+1) = By(t), \quad (x(0),y(0)) = (x_0, y_0)
\]
defined by nonnegative $n \times n$ matrices $A$, $B$.  Each system could describe the \emph{regional} dynamics of a structured population on a patch or location (for example).  We assume that the origin is an asymptotically stable equilibrium for both $\Sigma_A$, $\Sigma_B$: equivalently, both $A$ and $B$ are Schur-stable.  We also assume that a set $\mathcal{D}$ of nonnegative matrices is specified where for all $D$ in $\mathcal{D}$, both $A-D$, $B-D$ are entry-wise nonnegative.  The elements of $\mathcal{D}$ describe the admissible diffusion matrices modelling the coupling of the systems: in the context of population dynamics, these matrices describe the movement of species between patches.  The problem is to determine conditions that guarantee that the overall $2n$-dimensional coupled system, 
\begin{eqnarray}\label{eq:sysCoup}
\Sigma_{A,B}:x(t+1) &=& (A-D)x(t) + D y(t) \\
\nonumber y(t+1) &=& D x(t) + (B-D)y(t), 
\end{eqnarray}
has an asymptotically stable equilibrium at the origin for all $D$ in $\mathcal{D}$.  This is equivalent to the matrix
\begin{equation}\label{eq:M}
M = \left(\begin{array}{c c}
A-D & D\\
D & B-D
\end{array}\right)
\end{equation}
satisfying $\rho(M) < 1$ for all $D \in \mathcal{D}$.  When this happens, we say that the coupled system $\Sigma_{A,B}$ is robustly diffusively stable (RDS) with respect to $\mathcal{D}$.  

The set $\mathcal{D}$ of admissible coupling matrices can vary depending on the context or application.  In \cite{deLee}, $\mathcal{D}$ consists of matrices that act diffusively on the invariant cone of the system. For the cone $\mathbb{R}^n_+$ considered here, these are given by diagonal matrices with nonnegative entries.  In \cite{Guiv}, diffusion or movement between patches was described by general nonnegative matrices.  In both \cite{deLee, Guiv} \emph{common linear copositive Lyapunov functions} \cite{MasSho} are used to establish conditions for RDS.  It was also noted in these references that general common quadratic Lyapunov functions are not sufficient for robust diffusive stability.  This naturally suggests the following question.  Are there other classes of common Lyapunov function whose existence is sufficient for robust diffusive stability?  We address this question in Sections \ref{sec:2} and \ref{sec:jlf} and show that for diagonal coupling/diffusion matrices common diagonal Lyapunov functions or joint linear copositive Lyapunov functions are each sufficient for robust diffusive stability.



\section{Robust diffusive stability and common diagonal Lyapunov functions}\label{sec:2}
Throughout this section, we consider a pair of positive LTI systems, $\Sigma_A$, $\Sigma_B$, where $A$ and $B$ are Schur-stable nonnegative matrices in $\mathbb{R}^{n \times n}_+$.  The set of admissible coupling matrices is given by $\mathcal{D}_{A,B}$, the set of nonnegative diagonal matrices $D \in \mathbb{R}^{n \times n}_+$ such that $A-D \geq 0$, $B-D \geq 0$.     

The results presented in \cite{Guiv} imply that $\Sigma_{A,B}$ is RDS with respect to $\mathcal{D}_{A,B}$ if there exists $v \gg 0$ with $v^T A \ll v^T$, $v^T B \ll v^T$; a corresponding result for more general continuous-time systems was proven in \cite{deLee}.  $V(x) = v^Tx$ then defines a common linear copositive Lyapunov function for $\Sigma_A$, $\Sigma_B$.  

In the next result, we show that if there exists a common diagonal solution to the inequalities given in 4 (or 5) of Theorem \ref{thm:SchStab} defined by matrices $A$ and $B$, $\Sigma_{A,B}$ is RDS with respect to $\mathcal{D}_{A,B}$. 

\begin{thm}\label{thm:diag}
Let $A, B$ in $\mathbb{R}^{n \times n}_+$ be Schur-stable matrices.  Consider the following three statements. 
\begin{itemize}
\item[(i)] There exists a diagonal matrix $E \succ 0$ with \begin{equation}\label{eq:comSte}A^TEA - E \prec 0, \quad B^TEB-E \prec 0.\end{equation}
\item[(ii)] There exists a diagonal matrix $E \succ 0$ with \begin{equation}\label{eq:comLya} (A-I)^T E + E(A-I) \prec 0, \quad (B-I)^T E + E(B-I) \prec 0.\end{equation} 
\item[(iii)] For any $D \in \mathcal{D}_{A,B}$, the matrix $M$ in \eqref{eq:M} is Schur-stable.
\end{itemize}
Then (i) $\implies$ (ii) $\implies$ (iii).
\end{thm}
\begin{proof}
We first show that (i) implies (ii).  Suppose $E \succ 0$ is a diagonal matrix satisfying \eqref{eq:comSte}.  By adding and subtracting $EA$ to $A^TEA-E$, we see that
\begin{eqnarray} \nonumber A^TEA-E &=& (A-I)^TEA + EA - E \\ &=& \label{ca1} (A-I)^TEA+E(A-I).\end{eqnarray} If we now add and subtract $(A-I)^TE$, we see that
\begin{eqnarray} \label{ca2}
(A-I)^TEA+E(A-I)=(A-I)^TE+E(A-I)+(A-I)^TE(A-I).
\end{eqnarray}
Putting everything together, it follows that
\begin{equation}\label{eq:StLy}A^TEA-E = (A-I)^TE+E(A-I)+(A-I)^TE(A-I).\end{equation}
As $A$ is Schur, $A-I$ is invertible. Therefore, as $E \succ 0$, $(A-I)^TE(A-I) \succ 0$, and it follows from \eqref{eq:StLy} that \[ (A-I)^TE+E(A-I) \prec A^TEA-E \prec 0.\] An identical argument shows that $(B-I)^TE+E(B-I) \prec 0$.  Thus (i) implies (ii) as claimed. 

We next show that (ii) implies (iii).  Let $E \succ 0$ be diagonal with $Q_A=(A-I)^T E + E(A-I) \prec 0$, $Q_B=(B-I)^T E + E(B-I) \prec 0$, and let $D \in \mathcal{D}_{A, B}$ be given.  Consider the diagonal matrix 
\[P=\left(\begin{array}{c c}
E & 0\\
0 & E
\end{array}\right)\]
and compute $Q_M=(M-I)^TP+P(M-I)$ where $M$ is given by \eqref{eq:M}.  A direct calculation shows that 
\begin{eqnarray*} \label{P1}
Q_M &=& \left(\begin{array}{cc}
	    (A-I-D)^T E + E(A-I-D) & 2ED \\
    	2ED & (B-I-D)^T E + E(B-I-D)
    \end{array}\right)  \\
&=& \left(\begin{array}{cc}
	    Q_A - 2ED & 2ED \\
    	2ED & Q_B - 2ED
    \end{array}\right).
\end{eqnarray*} 
$Q_M$ is the sum of the negative definite, block-diagonal matrix, $\textrm{diag}(Q_A, Q_B)$, and the matrix 
 \[ \left(\begin{array}{cc}
	    - 2ED & 2ED \\
    	2ED & - 2ED
    \end{array}\right)\]
    which is easily seen to be negative semi-definite by an application of Gershgorin's theorem \cite{HorJoh}.  Thus $Q_M \prec 0$, and hence $M$ is Schur-stable by condition 5 of Theorem \ref{thm:SchStab}.  This completes the proof. 
\end{proof}

Note that the implication (i) $\implies$ (iii) means that the existence of a CDLF for $\Sigma_A$, $\Sigma_B$ is sufficient for RDS with respect to $\mathcal{D}_{A, B}$.  We next present a numerical example to show that conditions (i) and (ii) in Theorem \ref{thm:diag} are not equivalent and to highlight that the less conservative condition (ii) is distinct to the conditions for RDS based on common linear copositive Lyapunov functions given in \cite{deLee, Guiv}.

\begin{example}\label{diagex} Consider \begin{eqnarray*} 
A =  \left(\begin{array}{cc}
	    0.1 & 1 \\
    	0 & 0
    \end{array}\right), \;
    B= \left(\begin{array}{cc}
	    0.1 & 0 \\
    	1 & 0
    \end{array}\right)\end{eqnarray*} 
Clearly $A$ and $B$ are Schur-stable as both have eigenvalues $0, 0.1$.  Choosing $E=I$, a straightforward calculation shows that $(A-I)^T + (A-I)$ and $(B-I)^T + (B-I)$ are both equal to \begin{equation}\left(\begin{array}{cc}
	    -1.8 &  1 \\
    1  & -2
    \end{array}\right) \prec 0.
\end{equation}
Thus condition (ii) of Theorem \ref{thm:diag} is satisfied and we can conclude that the system $\Sigma_{A,B}$ is robustly diffusively stable with respect to $\mathcal{D}_{A,B}$. 

We now show by contradiction that $A$ and $B$ do not have a common linear copositive Lyapunov function.  Suppose that there exists $v \gg 0$ with $v^TA \ll v^T$, $v^TB \ll v^T$.  The first inequality means that $(0.1 v_1 \; v_1) \ll (v_1 \; v_2)$ implying that $v_1 < v_2$.  The second is equivalent to $(0.1v_1 + v_2 \; 0) \ll (v_1 \; v_2)$ implying that $v_1 > v_2$, which is a contradiction. This shows that no such $v$ exists.  So (ii) in Theorem \ref{thm:diag} gives a distinct condition for RDS to that based on the existence of a common linear copositive Lyapunov function.  

Finally, we also note that while there exists a diagonal $E \succ 0$ satisfying \eqref{eq:comLya}, there is no solution $E\succ 0$ to \eqref{eq:comSte}.  To see this, let \[E = \left(\begin{array}{c c}
1 & 0\\
0 & \varepsilon
\end{array}\right)
\] with $\varepsilon > 0$.  Then
\[A^TEA-E = \left(\begin{array}{c c}
-0.99 & 0.1\\
0.1 & 1-\varepsilon
\end{array}\right) ,
B^TEB-E = \left(\begin{array}{c c}
\varepsilon-0.99 & 0\\
0 & -\varepsilon
\end{array}\right) .
\]
If $A^TEA-E \prec 0$, we must have $\varepsilon > 1$ which implies that $B^TEB-E$ has a positive eigenvalue $\varepsilon-0.99$.  This shows that there is no common diagonal solution of \eqref{eq:comSte}.  Thus, condition (ii) does not imply condition (i) and (ii) is a strictly less conservative condition for RDS than (i).  
   \end{example}

\section{Joint linear copositive Lyapunov functions and RDS} \label{sec:jlf}
In this section, we again consider a pair of positive LTI systems $\Sigma_A$, $\Sigma_B$ whose system matrices $A$, $B$ are Schur-stable, and the set of coupling matrices  given by $\mathcal{D}_{A,B}$.  Joint linear copositive Lyapunov functions (JLCLFs) were introduced in \cite{WuSun} for the stability analysis of switched positive linear systems, with the aim of obtaining less conservative conditions for stability than those based on CLCLFs.  We show in Theorem \ref{thm:jlf} that JLCLF existence also gives a sufficient condition for RDS when at least one of $A$, $B$ is irreducible.  First, we recall the appropriate definition of a JLCLF in the discrete-time setting considered here.  

\begin{defn}\label{def:jlf}
If there exists a vector $v\gg 0$ with 
\begin{enumerate}
\item $v^T A \leq v^T$, $v^TB \leq v^T$;
\item $v^T(A+B) \ll 2v^T$,
\end{enumerate}      
then $V(x) = v^Tx$ defines a JLCLF for the systems $\Sigma_A$, $\Sigma_B$.  
\end{defn}
 
\begin{thm} \label{thm:jlf} Let $A$, $B$ be Schur-stable matrices in $\mathbb{R}^{n \times n}_+$, and assume that at least one of $A$, $B$ is irreducible.  Suppose that there exists a JLCLF for $\Sigma_A$, $\Sigma_B$.  Then $\Sigma_{A,B}$ is RDS with respect to $\mathcal{D}_{A,B}$. 
\end{thm}
\begin{proof}
Let $V(x)=v^Tx$ be a JLCLF for $\Sigma_A$, $\Sigma_B$.  We will argue by contradiction, so assume that there is some $D \in \mathcal{D}_{A,B}$ such that the matrix $M$ in \eqref{eq:M} has $\rho(M) = 1$.  By the Perron-Frobenius Theorem, there must exist some $y := (u \ w)^T > 0$, where $u, w$ are in $\mathbb{R}^n_+$, with $My=y$.  In terms of $u, w$, this means that
    \begin{eqnarray*}
    Au - Du +Dw &=& u\\
    Du + B w - D w &=& w.
    \end{eqnarray*}
If $w=0$, the second equation implies that $Du=0$ and the first equation then implies that $Au=u$ which contradicts the assumption that $A$ is Schur-stable.  Thus $w \neq 0$.  Similarly, we can show that $u \neq 0$.  A simple rearrangement of the above equations shows that 
\begin{eqnarray}\label{eq:AB1}
(A-I)u &=& D(u-w)\\
\nonumber (B-I)w &=& D(w-u).
\end{eqnarray}
Adding these shows that $(A-I)u+(B-I)w = 0$, which implies that 
\begin{equation}\label{eq:vAB}
v^T(A-I)u + v^T(B-I)w = 0.
\end{equation}
Now, note that from the definition of JLCLF, it follows that $v^T(A-I) \leq 0$, $v^T(B-I) \leq 0$ and $v^T(A-I)+v^T(B-I) \ll 0$.   Recall that $u \geq 0$, $w \geq 0$.  Suppose that  the $i^{th}$ component of $u$ is positive, $u_i > 0$.  Then as $v^T(A-I) \leq 0$, \eqref{eq:vAB} implies that the $i^{th}$ component of $v^T(A-I)$ must equal 0.  $v^T(A-I)+v^T(B-I) \ll 0$ now implies that the $i^{th}$ component of $v^T(B-I)$ is negative.  Finally, \eqref{eq:vAB} implies that $w_i = 0$.  A similar argument starting from $w$ establishes the following.  
For all $i$, $1\leq i \leq n$, at least one of $u_i$, $w_i$ must equal 0.  It follows that we can partition $\{1, \ldots, n\}$ into three (disjoint) sets $I_u$, $I_w$, $I_0$, where 
\begin{itemize}
\item $u_i > 0$, $w_i = 0$ for all $i \in I_u$;
\item $u_i=0$, $w_i > 0$ for all $i \in I_w$;
\item $u_i=w_i=0$ for all $i \in I_0$. 
\end{itemize}
As $u$ and $w$ are non-zero, $I_u$ and $I_w$ are both non-empty.  Choose any $i \in I_u$.  The second equation in \eqref{eq:AB1} implies that.
\[\sum_{j \in I_w} b_{ij} w_j = -d_i u_i.\]
As $w_j > 0$ for $j \in I_w$, and $d_iu_i \geq 0$, it follows that $b_{ij} = 0$ for all $j\in I_w$.  Next note that for all $i \in I_0$,
\[\sum_{j \in I_w} b_{ij} w_j = 0,\]
so again we can conclude that $b_{ij} = 0$ for all $j \in I_w$.  Putting this together, we see that $b_{ij} = 0$ for all $i \in I_u \cup I_u$, $j \in I_w$, which implies that $B$ is reducible.  By considering the first equation in \eqref{eq:AB1} instead, we can show in an identical fashion that $A$ is reducible. Thus, we have shown that if there exists $D \in \mathcal{D}_{A,B}$ such that $M$ given by \eqref{eq:M} is not Schur-stable, then both $A$ and $B$ are reducible.  So if one of $A$ or $B$ is irreducible, $M$ is Schur-stable for all $D \in \mathcal{D}_{A, B}$ and the result is proven.    

\end{proof}

The next example illustrates that systems $\Sigma_A$, $\Sigma_B$ can have a JLCLF without having a CLCLF.  This shows that Theorem \ref{thm:jlf} gives a less conservative condition for RDS (in the case where at least one matrix is irreducible). 

\begin{example}
Consider the matrices
\[A=\left(\begin{array}{c c}
0 & \frac{3}{4}\\
1 & 0
\end{array}\right), \;\; 
B =\left(\begin{array}{c c}
0 & 1\\
\frac{3}{4} & 0
\end{array}\right) .
\]
It is not hard to see that there is no $v \gg 0$ with $v^T A \ll v^T$, $v^T B \ll v^T$, as  $v^T A \ll v^T$ implies $v_2 < v_1$ while  $v^T B \ll v^T$ implies $v_1 < v_2$.  On the other hand, $v^T = (1 \;\; 1)$ clearly satisfies  $v^T A \leq v^T$,  $v^T B \leq v^T$ and $v^T(A+B) \leq 2v^T$, so Theorem \ref{thm:jlf} allows us to conclude that $\Sigma_{A,B}$ is RDS with respect to $\mathcal{D}_{A,B}$.
\end{example}

\section{Leslie matrices and RDS} \label{sec:Les}
In this section, we consider coupled LTI systems $\Sigma_A$, $\Sigma_B$ whose system matrices are \emph{extended Leslie matrices} \cite{Cus}.  Formally, we assume that $A$ and $B$ are nonnegative matrices of the form

\begin{eqnarray}\label{eq:LesAB}
\scalemath{0.95}{    A= \left( 
        \begin{array}{ccccc}
           a_{11} & a_{12} & a_{13} & \cdots & a_{1n} \\
		    a_{21} & 0 & 0 & \cdots & 0 \\
			0 & a_{32} & 0 & \cdots & 0 \\
			\vdots &  & \ddots  & \ddots & \vdots \\
			0 & 0 & 0 & a_{n n-1} & a_{nn}
		\end{array} \right), B= \left( 
        \begin{array}{ccccc}
            b_{11} & b_{12} & b_{13} & \cdots & b_{1n} \\
		    b_{21} & 0 & 0 & \cdots & 0 \\
			0 & b_{32} & 0 & \cdots & 0 \\
			\vdots &  & \ddots  & \ddots & \vdots \\
			0 & 0 & 0 & b_{n n-1} & b_{nn}
		\end{array} \right).}
		\end{eqnarray}
Leslie matrices of this form are often used in age-structured population models.  The number of age classes is $n$, and the $i^{th}$ component of the state vector, $x_i$, describes the size of the $i^{th}$ age class.  In such models, it is usually assumed that the time-step is equal to the length of time spent in each age class before moving to the next.  If the $(n,n)$ entry of the matrix is non-zero, a proportion of the oldest class survive to the next census.  The structure of the model means that diagonal matrices are not appropriate for modelling movement or diffusion for these models as if a member of the $i^{th}$ age cohort at time $t$ moves, then they will join the $i+1^{st}$ age cohort on the other patch at time $t+1$.  The only exception to this is the oldest cohort, where individuals who move patch would remain in the oldest cohort.  

Mathematically, this means that nonnegative matrices of the form
\begin{eqnarray} \label{eq:LesD}  D = \left( 
        \begin{array}{ccccc}
           0 & 0 & 0 & \cdots & 0 \\
		    d_{21} & 0 & 0 & \cdots & 0 \\
			0 & d_{32} & 0 & \cdots & 0 \\
			\vdots &  & \ddots  & \ddots & \vdots \\
			0 & 0 & 0 & d_{n n-1} & d_{nn}
		\end{array} \right), \end{eqnarray}
describe movement or dispersion between patches.  For the remainder of this section, we assume that $A$ and $B$ are of the form \eqref{eq:LesAB}, and use $\mathcal{L}_{A, B}$ to denote the set of all nonnegative matrices $D$ of the form \eqref{eq:LesD} such that $A-D \geq 0$, $B-D \geq 0$.  In the next result, we consider coupling/diffusion matrices belonging to the subset $\mathcal{L}^{1}_{A,B}$ of $\mathcal{L}_{A, B}$ consisting of matrices $D \in \mathcal{L}_{A, B}$ containing only one non-zero row.  From an applied perspective, this corresponds to one of the following situations.
\begin{enumerate}
\item Movement/diffusion only happens from one age class, such as is the case for some insect species such as moths, beetles, or butterflies.
\item Movement only happens from the two oldest age classes.
\end{enumerate}

\begin{thm} \label{thm:Les} Let $A, B  \in \mathbb{R}^{n \times n}_+$ be Schur-stable, Leslie matrices of the form \eqref{eq:LesAB}. Then the coupled system $\Sigma_{A,B}$ is RDS with respect to $\mathcal{L}^{1}_{A,B}$.
\end{thm}

\begin{proof}
        
We prove our result by contradiction.  Assume that $\Sigma_{A,B}$ is not RDS with respect to $\mathcal{L}^1_{A,B}$.  As $A$ and $B$ are Schur-stable, and the spectral radius is a continuous function of the entries of a matrix, it follows that there must exist some $D \in \mathcal{L}^{1}_{A,B}$ such that $\rho(M) = 1$ with $M$ given by \eqref{eq:M}.  From the Perron-Frobenius theorem, there exists some $y := (v \ w)^T > 0$, where $v, w$ are in $\mathbb{R}^n_+$, with $My=y$.  In terms of $v, w$, this means that
    \begin{eqnarray*}
    Av - Dv +Dw &=& v\\
    Dv + B w - D w &=& w.
    \end{eqnarray*}
    If $w=0$, then $Dv=w=0$ and the first equation implies that $Av=v$ which contradicts that $A$ is Schur-stable.  Thus $w \neq 0$.  Similarly, we can see that $v \neq 0$.  Next rearrange the above equations slightly as follows.
     \begin{eqnarray}
        \label{rhs1} &&(A-I)v = D(v-w)\\
        \label{rhs2} &&(B-I)w= D(w-v) .
    \end{eqnarray} As $D$ has only one non-zero row, it follows that either $D(v-w) \geq 0$ or $D(w-v) \geq 0$. In turn, this means that either $Av \geq v$ or $Bw \geq w$, which contradicts the assumption that both $A$ and $B$ are Schur-stable, as we have already shown that $v> 0$, $w > 0$. This contradiction shows that $M$ is Schur for all $D \in\mathcal{L}^{1}_{A,B}$ as claimed.
\end{proof}

In the light of the previous result, it is natural to ask whether robust diffusive stability also holds when movement is possible from two or more age classes, so that the matrix $D$ describing dispersal has two or more non-zero rows.  The following numerical example, generated using MATLAB \cite{Matlab} shows that this is not the case.  
\begin{example}
    Consider \begin{eqnarray*}
    A = \left(\begin{array}{ccc}
    0.1  & 0.85  &  0.15\\
    0.9  &    0  &    0\\
         0  & 0.7  &  0.2
    \end{array}\right), B = \left(\begin{array}{ccc}
 0.6 &  0.1  &  1\\
    0.5 &  0 &   0 \\
     0  & 0.35  &  0.45
         \end{array}\right).
    \end{eqnarray*} Then, it can be checked numerically that both $A$ and $B$ are Schur. Letting 
\begin{eqnarray*} D = \left(\begin{array}{ccc}
         0   &      0   &      0 \\
    0.25    &     0  &       0 \\
         0  &  0.2   &      0
         \end{array}\right)
    \end{eqnarray*} one can verify that $\rho(M) \approx 1.02 >1$.  Thus, when movement can take place from more than one class, RDS is no longer guaranteed.  This example demonstrates what is referred to as \emph{dispersal driven growth} in \cite{Guiv}.  For more on dispersal driven growth and related topics, see \cite{BlakeThesis}.
\end{example}

While Theorem \ref{thm:Les} only guarantees RDS with respect to $\mathcal{L}^1_{A, B}$, Proposition \ref{pro:Les2} below gives a simple sufficient condition for RDS with respect to the larger set $\mathcal{L}_{A,B}$.  

Given two matrices $A$, $B$ in $\mathbb{R}^{n \times n}_+$, a row selection is any matrix $C$ in $\mathbb{R}^{n \times n}_+$ where for $1 \leq i \leq n$, the $i^{th}$ row in $C$ is either the $i^{th}$ row of $A$ or the $i^{th}$ row of $B$.  Let $R(A,B)$ denote the set of all row selections formed from $A$ and $B$.  The next result is an immediate consequence of results in \cite{MasSho, Seet}. 
\begin{thm} \label{thm:rowsel}
Let $A, B$ be Schur-stable matrices in $\mathbb{R}^{n \times n}_+$.  There exists $v \gg 0$ in $\mathbb{R}^n_+$ with $Av \ll v$, $Bv \ll v$ if and only if every $C$ in $R(A,B)$ is Schur-stable.
\end{thm}
In the next result, given a pair $A,B $ of Leslie matrices of the form \eqref{eq:LesAB} the matrices $S_1$, $S_2$ are given by
\begin{eqnarray*}\scalemath{0.95}{
    S_1 = \left( 
        \begin{array}{ccccc}
           a_{11} & a_{12} & a_{13} & \cdots & a_{1n} \\
		    m_{21} & 0 & 0 & \cdots & 0 \\
			0 & m_{32} & 0 & \cdots & 0 \\
			\vdots &  & \ddots  & \ddots & \vdots \\
			0 & 0 & 0 & m_{n n-1} & m_{nn}
		\end{array} \right), S_2 = \left( 
        \begin{array}{ccccc}
            b_{11} & b_{12} & b_{13} & \cdots & b_{1n} \\
		    m_{21} & 0 & 0 & \cdots & 0 \\
			0 & m_{32} & 0 & \cdots & 0 \\
			\vdots &  & \ddots  & \ddots & \vdots \\
			0 & 0 & 0 & m_{n n-1} & m_{nn}
		\end{array} \right)}
		\end{eqnarray*}
where $m_{ii-1} := \max\left\{a_{ii-1}, b_{ii-1}\right\}$, for $ i \in \{2, ..., n\}$, and $m_{nn} :=  \max\left\{a_{nn}, b_{nn}\right\}$.

We now state a simple sufficient condition for the system $\Sigma_{A,B}$ to be RDS with respect to $\mathcal{L}_{A,B}$, when $A$, $B$ are Leslie matrices of the form \eqref{eq:LesAB}.  

\begin{pro}\label{pro:Les2}
Let $A, B \in \mathbb{R}^{n \times n}_+$ be Leslie matrices of form \eqref{eq:LesAB}.  Assume that $S_1, S_2$ are both Schur-stable. Then $\Sigma_{A,B}$ is RDS with respect to $\mathcal{L}_{A, B}$. 
\end{pro}
\begin{proof}
The proof is a simple application of Theorem \ref{thm:rowsel}.  From the form of $S_1, S_2$, any $C$ in $R(A,B)$ satisfies either $C \leq S_1$ or $C \leq S_2$.  As $C \geq 0$ and both $S_1, S_2$ are Schur-stable, this implies that $C$ is Schur-stable.  Theorem \ref{thm:rowsel} now implies that there exists $v \gg 0$ with $Av \ll v$, $Bv \ll v$.  Writing $\hat{v} = (v^T \;\; v^T)^T$, it is easy to see that for all $D \in \mathcal{L}_{A,B}$, $M \hat v \ll \hat{v}$ where $M$ is given by \eqref{eq:M}.  Thus $M$ is Schur-stable by the remark after Theorem \ref{thm:SchStab}.
\end{proof} 
  
\section{Concluding remarks}
Building on the work recently published in \cite{deLee,Guiv}, we have presented several novel results on robust diffusive stability (RDS) for coupled discrete-time positive LTI systems.  In particular, Theorem \ref{thm:diag} shows that common diagonal Lyapunov functions can be used to establish robust diffusive stability (RDS) when coupling matrices are diagonal.  Theorem \ref{thm:jlf} extends the conditions for RDS based on common linear copositive Lyapunov functions by showing that the weaker condition of joint linear copositive Lyapunov function existence is sufficient for RDS when at least one of the system matrices is irreducible.  The results of Section \ref{sec:Les} describe conditions for RDS for systems arising in population dynamics when the regional systems are described by extended Leslie matrices.  Theorem \ref{thm:Les}, in particular, shows that movement cannot lead to instability when it only occurs from one age class.  In future, it would be interesting to investigate extensions to nonlinear models, and to determine necessary conditions for RDS, which could help identify sufficient conditions for dispersal-driven growth \cite{Guiv}. 
\section*{Declaration of competing interest}
There are no competing interests to declare.
\section*{Acknowledgements}
BMC was supported by a Taighde Ireland postgraduate scholarship (GOIPG/2020/939)







\end{document}